\documentclass[11pt]{article}
\usepackage{amssymb, amsmath, amsthm, amsfonts,xcolor, enumerate}
\usepackage{fullpage}
\usepackage{hyperref}

\usepackage{t1enc}
\usepackage[T1]{fontenc}
\usepackage[utf8]{inputenc}

\newtheorem{theorem}{Theorem}[section]
\newtheorem{lemma}[theorem]{Lemma}
\newtheorem{obs}[theorem]{Observation}
\newtheorem{defn}[theorem]{Definition}
\newtheorem{cor}[theorem]{Corollary}

\newtheorem{claim}[theorem]{Claim}
\newtheorem{fact}[theorem]{Fact}

\newcommand{\ep}{\varepsilon}

\def\qedf{\hfill $\Box$}

\def\komment#1{}
\let\komment=\footnote

\begin{document}

\title{On the Ramsey-Turán problem for 4-cliques }

\author{ B\'ela Csaba\thanks{E-mail: bcsaba@math.u-szeged.hu. This research was supported
by the Ministry of Innovation and Technology of Hungary from the National Research, Development and
Innovation Fund, project no. TKP2021-NVA-09.}\\ Bolyai Institute, University of Szeged, Hungary}

\date{}

\maketitle

\begin{abstract}
We present an essentially tight bound for the Ramsey-Tur\'an problem for 4-cliques without using the Regularity lemma. 
This enables us to substantially extend the range in which one has the tight bound for the number of edges in $K_4$-free graphs as a function of 
the independence number, apart from lower order terms. 
\end{abstract}

\section{Introduction}

Two classical results in combinatorics are the Ramsey theorem~\cite{Ry} and the Turán~\cite{Tu} theorem. In 1969 Vera T.~Sós in~\cite{TS} introduced a common
generalization of these, the so called Ramsey-Turán problems. 
The Ramsey-Tur\'an number $RT(n,H,m)$
is the maximum number of edges a graph $G$ on $n$ vertices with independence
number less than $m$ can have without containing $H$ as a subgraph. In this paper we consider the $H=K_4$ case. For more information on Ramsey-Turán theory the reader may consult with
the excellent survey~\cite{STS} by Simonovits and Sós.

It is convenient to introduce the notation $\alpha=m/n,$ we will use it throughout the paper. In 1972 Endre Szemerédi~\cite{SzRT} proved the following upper bound for $K_4$-free graphs
with a small independence number.

\begin{theorem}\label{elso} For every $\eta>0$ there exists an $\alpha>0$ such that the following holds. Every $n$-vertex graph
having at least $(\frac{1}{8}+\eta) n^2$ edges contains either a $K_4$ or an independent set larger than $\alpha n.$ 
\end{theorem}

This result turned out to be almost tight, Bollobás and Erdős~\cite{BE} constructed  $K_4$-free graphs with independence number $o(n)$ having $n^2/8 -o(n^2)$ edges. 

The proof\footnote{This 
was the first occurrence of a (weak) regularity notion.} of Theorem~\ref{elso} gives the upper bound of $O((\log \log \frac{1}{\alpha})^{-1/2 + o(1)})$ for $\eta.$  Conlon and Schacht, independently, in unpublished work, showed that the bound
$O((\log \frac{1}{\alpha})^{-1/2 + o(1)})$ is already sufficiently large for $\eta,$ they applied the Frieze-Kannan weak regularity 
lemma~\cite{FriezeKannan}. 

In~\cite{FLZ}  Fox, Loh and Zhao gave significantly better upper and lower bounds for $\eta$ as a function of $\alpha.$ Their paper contains
several deep results concerning the Ramsey-Turán problem for $4$-cliques, here we only mention those which are most relevant for our main theorem. 
Fox et al.~proved that there exists
an absolute constant $\gamma>0$ such that if $\alpha \le \gamma,$ $\alpha(G)=\alpha n,$ and $e(G)\ge \frac{n^2}{8}+ \frac{3}{2}\alpha n^2,$ then $G$ has a $K_4$ (Theorem 1.6 in~\cite{FLZ}). 
  
For the lower bounds in~\cite{FLZ} set $\beta=\sqrt{(\log\log n)^3/\log n}$ and assume that $0<\alpha <1/3$ such that $\alpha/\beta \longrightarrow \infty.$ 
They constructed (Theorem 1.7) $K_4$-free 
graphs with maximum independent set of size $\alpha n$ and number of edges at least $\frac{n^2}{8}+(\frac{1}{3}-o(1))\alpha n^2.$ Moreover, in the case $\alpha^2/\beta\longrightarrow \infty$ 
they gave an even stronger lower bound for the number of edges: $\frac{n^2}{8}+(\alpha-\alpha^2)\frac{n^2}{2}-\beta n^2.$
That is, they proved the following bounds for $\eta$ if $\alpha$ is sufficiently larger than $\beta$:
$$\frac{1}{2}(\alpha - \alpha^2)-o(\alpha^2) \le \eta \le \frac{3}{2}\alpha.$$
We remark that for the upper bound they used the Regularity lemma~\cite{SzR}.

Recently, Lüders and Reiher~\cite{LR} proved\footnote{In fact they considered the Ramsey-Turán problem for every $K_r$ ($r\ge 3$), and proved sharp bounds.} that the above lower bound is essentially the truth in case $\alpha$ is a constant (although very small).

\begin{theorem}\label{LRe}
There exists a threshold $0<\gamma^*\ll 1$  
such that if $0<\gamma \le \gamma^*,$ $n$ is sufficiently large, $G$ is a graph on $n$ vertices with 
$e(G)>(n^2+n)/8 + (\gamma -\gamma^2)n^2/2$ and $\alpha(G)\le \gamma n,$ then $G$ contains a $K_4.$ 
\end{theorem}

Observe, that the above bound for $\eta$ is $(\gamma - \gamma^2)/2,$ which could be much larger than $(\alpha-\alpha^2)/2,$ whenever $\alpha \longrightarrow 0$ as $n \longrightarrow \infty.$
Since in the proof they used the Regularity lemma, $\gamma^*$ is very small, and the sufficiently large $n$ is in fact a tower-type function of $\gamma.$


There is a significant interest in finding new, ``Regularity lemma-free'' proofs for theorems that use the Regularity lemma. This motivated our investigations. In the present paper we improve upon Theorem~\ref{LRe} in two ways. First, the value we give for $\eta$ is essentially
optimal for any $\alpha$ which is smaller than an absolute constant. Second, since we do not use the Regularity lemma,
the constants are much better, only single-exponential. 
Our main result is the following:

\begin{theorem}\label{fottetel}
Set $\nu = 1/500,$ and let $\gamma=\exp(- 10 \log (1/\nu)/\nu)$ and 
$N= \exp(10\log (1/\nu)/\nu).$ If $n\ge N,$ $G$ is a graph on $n$ vertices with $\alpha=\alpha(G)/n \le \gamma$ and 
$$e(G)>\frac{n^2+n}{8} + (\alpha -\alpha^2)\frac{n^2}{2},$$ then $G$ contains a $K_4.$ 
\end{theorem}

We did not optimize for the value of $\nu$ in Theorem~\ref{fottetel}, it is possible that with careful
computation one can significantly improve upon it. 

Several ideas, which play an essential role in the proof of Theorem~\ref{fottetel}, are taken from the papers~\cite{FLZ} and~\cite{LR}, together with a
recent graph decomposition method of the author~\cite{Cs}, which has some common features with a key lemma in~\cite{SzRT}.

The outline of the paper is as follows. In the second section we discuss the necessary notions, and prove a key lemma. In the final section we prove our main theorem. 
Let us remark that we will not be concerned with floor signs, divisibility, and so on in ~the proof. This makes the notation simpler, and the proof
easier to follow.

\section{Definitions, main tools}\label{defek}

Given a graph $G$ with vertex set $V$ and edge set $E$ we let $e(G)=|E(G)|.$ For a vertex $v\in V$ the degree of $v$ is denoted by $deg_G(v),$ if it is clear from
the context,
the subscription may be omitted. The neighborhood of $v$ is denoted by $N(v),$ so $deg(v)=|N(v)|.$ The minimum degree of $G$ is denoted by $\delta(G).$ If $S\subset V,$ then
$G[S]$ denotes the subgraph of $G$ induced by $S$ and $deg(v, S)=|N(v)\cap S|.$ Given two disjoint
sets $S, T\subset V$ the bipartite subgraph induced by them is denoted by $G[S, T],$ and $E(G[S, T])=E(S, T),$ and we let $e(S, T)=|E(S, T)|.$

Suppose that $F=(V, E)$ is a graph with non-empty subsets $A, B\subset V,$ $A\cap B=\emptyset.$ Then the density of $F[A, B]$ is $$d_F(A, B)=\frac{e_F(A, B)}{|A|\cdot |B|}.$$

\subsection{On the minimum degree of $G$}\label{minifok}

Let $G=(V, E)$ be a graph on $n$ vertices with independence number $\alpha(G) =\alpha n,$ where $0<\alpha<1/3.$ 
We assume, as in Theorem~\ref{fottetel}, that $$e(G)> \frac{n^2+n}{8}+(\alpha -\alpha^2)\frac{n^2}{2}.$$

Take the smallest subset $V_1\subset V$ for which 
$$e(G[V_1])>\frac{|V_1|^2+|V_1|}{8}+(\alpha -\alpha^2)\frac{n^2}{2}$$ (note, that we have $n,$ not $|V_1|$ in the second part 
of the expression on the right). Let us denote $G[V_1]$ by $G_1.$ 

\begin{claim}\label{meret}
The number of vertices in $G_1$ is more than $\sqrt{\alpha(G) n/2}\ge \sqrt{\alpha(G) N/2}.$ 
\end{claim}

\noindent {\bf Proof:} It is easy to see, that $e(G_1)>  (\alpha -\alpha^2)\frac{n^2}{2}.$ Since $\alpha\le \gamma <1/2,$ we have $e(G_1)> {\alpha n^2}/4=\alpha(G)n/4.$ Using that $G_1$ is a simple graph, this implies that 
it must have more than $\sqrt{\alpha(G) n/2}$ vertices. The fact that $n\ge N$ finishes the proof. \qedf

\medskip

Set $n_1=|V_1|,$ and let $\alpha_1=\alpha(G_1)/n_1.$ Clearly, $\alpha(G_1)\le \alpha(G),$ as $G_1$ cannot have a larger independent set than what $G$ has.


\begin{claim}\label{minfok}
The minimum degree of $G_1$ is at least $n_1/4.$ Moreover, $$e(G_1)> \frac{n_1^2+n_1}{8}+ (\alpha_1 -\alpha_1^2)\frac{n_1^2}{2}.$$
\end{claim}

\noindent {\bf Proof:} 
Using that $V_1$ is the smallest subset with the above property, for 
every $v\in V_1$ we have $$e(G_1[V_1-v])\le \frac{(n_1-1)^2+n_1-1}{8}+(\alpha -\alpha^2)\frac{n^2}{2}.$$ 
Hence, $$deg_{G_1}(v)=e(G_1[V_1])-e(G_1[V_1-v])\ge  \frac{(n_1)^2+n_1}{8}+(\alpha -\alpha^2)\frac{n^2}{2} -
 \frac{(n_1)^2-n_1}{8}-(\alpha -\alpha^2)\frac{n^2}{2}=\frac{n_1}{4}.$$

Clearly, for the second part of the claim it is enough to prove that $$\alpha(G)n-\alpha(G)^2\ge \alpha(G_1)n -\alpha(G_1)^2,$$
since $n_1\le n$ implies $\alpha(G_1)n -\alpha(G_1)^2\ge  \alpha(G_1)n_1 -\alpha(G_1)^2.$
Dividing by $n^2$ we get $$\frac{\alpha(G)}{n}-\frac{\alpha(G)^2}{n^2}\ge \frac{\alpha(G_1)}{n}-\frac{\alpha(G_1)^2}{n^2}.$$
Since $\alpha(G_1)\le \alpha(G)<n/2$ and the function $f(x)=x-x^2$ is monotone increasing in the interval $(0, 1/2),$ the above inequality
is satisfied, proving what was desired.
\qedf

\medskip 

Observe, that we may suppose that $\alpha(G)\ge 2,$ otherwise $G$ is a complete graph, and the
theorem becomes trivial. 
Using this observation and the above claims, in proving Theorem~\ref{fottetel} we can restrict our attention to graphs of order $n\ge \sqrt{\alpha(G) N/2}\ge \sqrt{N}= \exp(5\log (1/\nu)/\nu)=\exp(2500 \cdot \log 500),$ having minimum degree at least $n/4.$ 

\subsection{Finding a large quasi-random pair}

Given a number $\varepsilon \in (0,1)$ we say that a bipartite graph $F$ with parts $A$ and $B$ is an {\it $\varepsilon$-regular pair} if the following holds for every 
$A'\subset A,$ $|A'|\ge \varepsilon |A|$ and $B'\subset B,$ $|B'|\ge \varepsilon |B|$: $$|d_F(A, B)-d_F(A', B')|\le \varepsilon.$$ 
This notion plays a central role in the Regularity lemma~\cite{SzR}. In this paper we work with a closely related
but more permissive, one-sided definition. 

\begin{defn}
Given a bipartite graph $F$ with parts $A$ and $B$ we say that $F$ is an {\em $\varepsilon^+$-regular pair}, if for any $A'\subset A, B'\subset B$ with $|A'|\ge \varepsilon |A|,$ $|B'|\ge \ep |B|$ we have $d_F(A', B')\ge \ep.$ 
\end{defn}

We need a simple fact which is called the {\it convexity of density} (see 
eg.~in~\cite{KS}), the proof is left for the reader.

\begin{claim}\label{convexity} Let $F=F(A, B)$ be a bipartite graph, and let $1\le k \le |A|$ and $1\le m \le |B|.$ Then 
$$d_F(A, B)=\frac{1}{\binom{|A|}{k}\binom{|B|}{m}}\sum_{X\in \binom{A}{k}, Y\in \binom{B}{m}} d(X, Y).$$
\end{claim}

In other words, ${\mathbb E}_{X, Y}d(X, Y)=d_F(A, B),$ where $X\subseteq A,$ $Y\subseteq B$ are randomly
chosen subsets with $|X|=k$ and $|Y|=m.$

\smallskip

In order to prove Theorem~\ref{fottetel} we need a lemma which plays a key role in the proof. This lemma asserts
that one can find a large bipartite quasi-random subgraph in a sufficiently dense graph, a much larger one than 
that given by the Regularity lemma. We remark that results of this kind were proved in~\cite{KS} using the graph functional method of Koml\'os, and in~\cite{Peng} by Peng, R\"odl and Ruci\'nski~\cite{Peng}.
Here the situation is different. One side of the quasi-random pair could be much bigger, the pair in general is
not balanced.

\begin{lemma}\label{egylepes}
Let $F$ be a bipartite graph  with parts $A$ and $B$ such that $|A|=a$ and $|B|=b,$ and every
vertex of $A$ has at least $\delta b$ neighbors in $B$ for some $\delta>0.$ Let $0<\varepsilon\le \delta/6$ be a real number. 
Then $F$ contains an $\varepsilon^+$-regular pair $F[X, Y]$ such that $|X|\ge  \exp\left(-2\log(\frac{2}{\varepsilon})\log(\frac{2}{\delta})/\ep \right)a$ and $|Y|\ge (\delta-2\ep)b.$ Furthermore, 
$deg(v, Y)\ge (\delta-2\ep)b$ for every $v\in X.$
\end{lemma}


\noindent {\bf Proof:} We prove the lemma by finding two sequences of sets $X_0, X_1, \ldots, X_l$ and $Y_0, Y_1,$ 
\ldots, $Y_l$ where $X_0=A,$ $Y_0=B$ and $F[X_l, Y_l]$ is $\varepsilon^+$-regular such that for every $1\le i\le l$ we have that $X_i\subset X_{i-1}$ and $Y_i\subset Y_{i-1}$ and 
$$\varepsilon |X_{i-1}|/2 \le |X_i|\le \varepsilon |X_{i-1}|$$ and $$|Y_i|=(1-\ep)|Y_{i-1}|.$$ 
Hence, we may choose $X=X_l$ and $Y=Y_l.$

We find the set sequences $\{X_i\}_{i\ge 1}$ and $\{Y_i\}_{i\ge 1}$ by the help of an iterative procedure. This procedure stops in the 
$i$th step, if $F[X_i, Y_i]$ is $\varepsilon^+$-regular. We have another stopping rule: if 
$|Y_i| \le (\delta (1+\ep/2)-2\ep)b$ for some $i,$ we stop. Later we will see
that in this case we have found what is desired, $F[X_i, Y_i]$ must be an $\varepsilon^+$-regular pair.

In the beginning we check if $F[X_0, Y_0]$ is an $\varepsilon^+$-regular pair. If it is, we stop. If not, then $X_0$ has a subset $S$ and $Y_0$ has a subset $T$ such that $|S|\ge \ep|X_0|,$ $|T|\ge \ep |Y_0|$ and $d_F(S, T)<\ep.$
Set $k=\ep|X_0|$ and $m=\ep|Y_0|,$ and apply Claim~\ref{convexity} for $F[S, T].$ We get that
the average density between $k$ element subsets of $S$ and $m$ element subsets of $T$ is less than
$\ep.$  Hence, there must exist subsets $X'_1\subseteq S$ with $|X'_1|=\varepsilon |X_0|$ and  
$Y'_1\subseteq T$ with $|Y'_1|=\ep|Y_0|$ such that $e(F[X'_1, Y'_1])<\ep |X'_1|\cdot |Y'_1|,$ 

We define a new set $X''_1\subset X'_1$ as follows: $$X''_1=\{x\in X'_1: \ deg(x, Y'_1) > 2\ep |Y'_1|\}.$$
Simple counting shows that $|X''_1|\le |X'_1|/2.$ Let $X_1=X'_1-X''_1,$ that is, those vertices of $X'_1$ that have at most
$2\ep |Y'_1|$ neighbors in $|Y'_1|.$ By the above we have $|X'_1|/2\le |X_1|\le |X'_1|.$ Set $Y_1=Y_0-Y'_1.$  
It is easy to see, that $\ep |X_0|/2\le |X_1|\le \ep|X_0|$ and $|Y_1|=(1-\ep)|Y_0|.$

For $i\ge 2$ the above is generalized. If $F[X_{i-1}, Y_{i-1}]$ is not  an $\varepsilon^+$-regular pair then we do the following. First, using Claim~\ref{convexity} we find $X'_i\subset X_{i-1}$ 
and $Y'_i\subset Y_{i-1}$ such that $|X'_i| =\varepsilon |X_{i-1}|$ and $|Y'_i| =\ep |Y_{i-1}|$ and $e(F[X'_i, Y'_i])<\ep |X'_i|\cdot |Y'_i|.$ Similarly to the above, we let  $$X''_i=\{x\in X'_i: \ deg(x, Y'_i) > 2\ep |Y'_i|\},$$ and conclude 
that $|X''_i|\le |X'_i|/2.$ Next we let $X_i=X'_i-X''_i,$ hence, $X_i$ is the set of those vertices of $X'_i$ that have at most $2\ep |Y'_i|$ neighbors in $Y'_i.$ 
Clearly, we have $|X'_i|/2\le |X_i|\le |X'_i|.$ Finally, we let $Y_i=Y_{i-1}-Y'_i.$ Hence, we have $\varepsilon |X_{i-1}|/2 \le |X_i|\le \varepsilon |X_{i-1}|$ and $|Y_i|=(1-\ep)|Y_{i-1}|.$
With this we proved that the claimed bounds for $|X_i|$ and $|Y_i|$ hold for every $i.$ 

\medskip

It might not be so clear that this process stops in a relatively few iteration steps.  
We need the following simple claim, the proof is left for the reader.

\begin{claim}\label{fok} 
Let $i\ge 1$ be an integer. If $u\in X_i,$ then $u$ has at most $2\ep (|Y'_1|+\ldots +|Y'_i|)$ 
neighbors in $B-Y_i.$ 
\end{claim}

Clearly, using the above claim we must have that $deg(v, Y_l)\ge (\delta - 2\ep)b$ for every $v\in X_l,$ which also implies that
$|Y_{l}|\ge (\delta -2\ep)b.$

Next we show that if $(\delta -2\ep)b\le |Y_l|\le (\delta(1+\ep/2)-2\ep)b$ then $F[X_l, Y_l]$ must be $\ep^+$regular. 
Assume that $u\in X_l.$ Then $deg(u, Y_l)\ge (\delta-2\ep)b,$ using Claim~\ref{fok}, hence, the number of 
{\it non-neighbors} of $u$ in $Y_l$ is at most $(\delta(1+\ep/2)-2\ep)b-(\delta -2\ep)b=\delta\ep b/2.$ 
Let $Y'\subset Y_l$ be arbitrary with $|Y'|=\ep |Y_l|$ and $u\in X_l.$ The number of neighbors of $u$ in $Y'$ is at least
$\ep |Y_l|-\delta \ep b/2.$ We will show that $$\ep |Y_l|-\delta \ep b/2\ge \ep |Y'|=\ep^2 |Y_l|.$$
This is equivalent to $$|Y_l|(1-\ep) \ge\delta\frac{ b}{2}.$$ Since $|Y_l|\ge (\delta -2\ep)b,$ it is sufficient if
$$(\delta-2\ep)(1-\ep) \ge \frac{\delta}{2}.$$ Using the condition $\ep \le \delta/6\le 1/6,$ (since $\delta\le 1$) we have
$$(\delta-2\ep) (1-\ep) \ge \frac{2\delta}{3 }\cdot \frac{5}{6}> \frac{\delta}{2},$$ as promised.  

Hence, for every $X'\subset X_l$ and $Y'\subset Y_l$ with $|Y'|=\ep |Y_l|$ we have $$e(X', Y')\ge \ep |X'|\cdot |Y'|.$$  
That is, if the procedure stopped because we applied the stopping rule, then the resulting pair must always be $\varepsilon^+$-regular. 

Next we upper bound the number of iteration steps. In every step the $Y$-side shrinks by a factor of $(1-\ep).$ 
We also have that $|Y_l|\ge (\delta-2\ep)b.$ Putting these together we get that
$$(1-\ep)^l\ge (\delta-2\ep)>\frac{\delta}{2}.$$ Hence, $$l< \frac{\log (2/\delta)}{\log(1/(1-\ep))}<\frac{\log(2/\delta)}{\ep},$$
here we used elementary calculus, in particular, that $e^{x}=\sum x^k/k!<\sum x^k =1/(1-x)$ for every real number $x\in (0, 1).$ 

Finally, we show the lower bound for $|X_l|.$ Note, that $|X_i|/|X_{i-1}|\ge \varepsilon/2$ for every $i\ge 1.$ Hence, 
$$|X_l| \ge \left(\frac{\varepsilon}{2}\right)^la =e^{-\log(2/\varepsilon)\log(2/\delta)/\ep}a.$$
 \qedf

\medskip

Let us remark that in~\cite{Cs} a similar lemma is proved for a more complicated definition of quasi-randomness.

\section{Proof of the main theorem}

Our main tool for proving Theorem~\ref{fottetel} is Lemma~\ref{egylepes}, but the following simple observation will also be useful. 

\begin{obs}\label{fokos}
Let $F=(V(F), E(F))$ be a $K_4$-free graph. Assume that $S\subset V(F)$ and $uv\in E(F)$ for some $u, v\in V(F).$
Then $deg_F(u, S)+deg_F(v, S)\le |S|+\alpha(F).$
\end{obs}

\medskip

Throughout we assume that $K_4 \not\subset G,$ and arrive at a contradiction, an upper bound for the number of edges in $G,$ which is less than
the lower bound for $e(G)$ in the theorem. 
We also assume that $\delta(G)\ge n/4,$ using the method in Section~\ref{minifok}. Recall, that we set
$\nu=1/500,$ and $\alpha(G)\le \gamma n=\exp(-5000 \cdot \log 500) \cdot n,$ where we substituted the value of $\nu$ in the expression for $\gamma.$ For convenience, we will keep 
the notation $\nu$ for the number $1/500.$

\smallskip

The fact below follows from the values of $\nu$ and $\gamma,$ we state it for future reference. 

\begin{fact}\label{kisalfa}
We have $\alpha = \alpha(G)/n<\exp(-\log (2/\nu)/\nu) \cdot \nu^2 < \nu^3.$
\end{fact} 

\medskip

We go through the proof of Theorem~\ref{fottetel} step-by-step, as follows. 
\medskip

\noindent {\bf Step 1.} Take an arbitrary subset $B_0\subset V$ such that $|B_0|=\nu n,$ set $A_0=  V-B_0,$ and apply Lemma~\ref{egylepes} for
the bipartite graph $G[B_0, A_0]$ with parameter $\nu.$ Note, that $deg(v, A_0)\ge n/4-\nu n$ for every $v\in B_0.$

We obtain a $\nu^+$-regular pair $G[B_1, A_1],$ where $B_1\subset B_0$ and 
$A_1\subset A_0.$ Fact~\ref{kisalfa} implies that the threshold numbers $\gamma$ and $N$ 
in the theorem are sufficiently large so that $|B_1|\ge \exp(-\log (2/\nu)/\nu)\nu n>\alpha(G),$ 
moreover, $d(v, A_1)\ge d(v, A_0)-2\nu n \ge n/4-3\nu n$ for every $v\in B_1.$
We also have the following.

\smallskip

\begin{claim}\label{A1nagy}
Since $G$ has no $K_4,$ we must have $|A_1|\ge n/2 -6\nu n-\alpha(G).$
\end{claim} 

\noindent {\bf Proof:} Since $|B_1|>\alpha(G),$ $B_1$ has at least one edge. Every vertex of $B_1$ has at least
$n/4-3\nu n$ neighbors in $A_1,$ hence, by Observation~\ref{fokos} we get the claimed bound for $|A_1|.$ \qedf

\medskip

\noindent {\bf Step 2.} We discard those vertices from $A_1$ that have less than $\nu |B_1|$ neighbors in $B_1$: let 
$$A_2=A_1-\{v: v\in A_1, \ deg(v, B_1)< \nu |B_1|\}.$$ By $\nu^+$-regularity we get that 
$|A_2|\ge (1-\nu)|A_1|\ge |A_1|-\nu n\ge n/2-7\nu n -\alpha(G),$ and $deg(v, A_2)\ge deg(v, A_1)-\nu n\ge  n/4-4\nu n$ for every $v\in B_1.$
We need the following.

\begin{claim}\label{A2fok}
For every vertex $v\in A_2$ we have $deg(v, A_2)\le \nu |A_1|.$   
\end{claim}

\noindent {\bf Proof:} Suppose not. Since $|N(v)\cap B_1| \ge \nu |B_1|,$ the density of edges between $N(v)\cap A_2$ and
$N(v)\cap B_1$ is at least $\nu$ using $\nu^+$-regularity. Hence, there exists a vertex $w\in N(v)\cap B_1$ having at least 
$\nu |N(v)\cap A_2|\ge  \nu^2 |A_1|>\nu^2n/3>\alpha(G)$ neighbors in $N(v)\cap A_2,$ where the last inequality follows from Fact~\ref{kisalfa}. But this contradicts with the 
$K_4$-freeness of $G,$ proving what was desired. 
\qedf

\medskip

\noindent {\bf Step 3.} 
Next we apply Lemma~\ref{egylepes} for the bipartite graph $G[A_2, V-A_2]$ with parameter $\nu.$ We have large minimum degree since by Claim~\ref{A2fok} every $v\in A_2$
has at least $n/4-\nu n$ neighbors in $V-A_2.$ 
We obtain a $\nu^+$-regular pair $(A_2', B_2),$ where $A_2'\subset A_2,$ $B_2\subset V-A_2,$ 
$|A_2'|\ge \exp(-2\log(2/\nu)/\nu)\nu n/3,$ and $|B_2|\ge n/2-6\nu n -\alpha(G),$ where the lower bound for the cardinality of $B_2$ follows from Observation~\ref{fokos} as in 
Claim~\ref{A1nagy}.

As before for $A_1,$ we discard those vertices of $B_2$ that have less than $\nu |A_2'|$ neighbors in $A_2'$:
 let 
$$B_3=B_2-\{v: v\in B_2, \ deg(v, A_2')< \nu |A_2'|\}.$$ By $\nu^+$-regularity we get that 
$|B_3|\ge (1-\nu)|B_2|\ge |B_2|-\nu n\ge n/2-7\nu n -\alpha(G).$ 
Using the arguments of Claim~\ref{A2fok} we get the following claim, the proof is very similar, we omit it.

\begin{claim}\label{B3fok}
For every vertex $v\in B_3$ we have $deg(v, B_3)\le \nu |B_2|.$   
\end{claim}

Summarizing, at this point we have a subgraph that is spanned by $A_2\cup B_3,$ inside $A_2$ or $B_3$ the vertices 
have only a few neighbors, and $A_2,$ $B_3$ have cardinality at least $n/2-7\nu n -\alpha(G).$ 

\medskip

\noindent {\bf Step 4.} Divide the set $V-(A_2\cup B_3)$ into two parts, $A'$ and $B',$ as follows:
$$A'=\{v: v\in V-(A_2\cup B_3), \ deg(v, B_3)\ge deg(v, A_2)\}$$ and $$ B'=\{v: v\in V-(A_2\cup B_3), \ deg(v, B_3)< deg(v, A_2)\}.$$

Since $|A_2|, |B_3|\ge n/2-7\nu n-\alpha(G),$ we have that $|A'|+|B'|\le 14\nu n+2\alpha(G)<15\nu n.$

\begin{obs}\label{atmenofok}
For every $v\in A'$ we have $deg(v, B_3)\ge n/9$ and $deg(u, A_2)\ge n/9$ whenever $u\in B'.$ This follows from the minimum
degree condition on $G$ and that $|A_2\cup B_3|\ge (1-15\nu)n.$
\end{obs}

\medskip

Let $A=A_2\cup A'$ and $B=B_3\cup B'.$ Clearly, $A\cap B=\emptyset,$ $A\cup B=V,$ and 
$$(1/2-7\nu)n-\alpha(G)\le |A|, |B|\le (1/2+7\nu)n+\alpha(G).$$
We have the following.

\begin{claim}\label{bentifok}
If $v\in A$ then $deg(v, A-A')\le \alpha(G).$ Similarly, $deg(u, B-B')\le \alpha(G),$ if $u\in B.$
\end{claim}

\noindent {\bf Proof:} We prove the first part of the statement, for vertices of $A,$ the second part can be proved similarly. Suppose on the 
contrary that $deg(v, A-A')>\alpha(G)$ for some $v\in A.$ Then there exists $u_1, u_2\in N(v)\cap (A-A')$ such that
$u_1u_2\in E.$ Since $G$ has no $K_4,$ $|N(u_1)\cap N(u_2)|\le \alpha(G).$ Both $u_1$ and $u_2$ can have up to
$|A'|\le 15\nu n$ neighbors in $A'$ and less than $\nu n$ neighbors in $A-A'.$ Hence, by Observation~\ref{fokos},
$$|(N(u_1)\cup N(u_2))\cap B| \ge n/2-32\nu n -\alpha(G)>n/2-33\nu n.$$
Putting these together we get
$$|B-(N(u_1)\cup N(u_2))|<(1/2 +7\nu)n+\alpha(G) - (1/2-33\nu)n<41\nu n.$$ Since every vertex of $A$ has at least
$n/9$ neighbors in $B,$ the above inequalities imply that either $|N(v)\cap N(u_1)\cap B|>\frac{1}{2}(n/9-41\nu n)>\alpha(G),$
or $|N(v)\cap N(u_2)\cap B|>\frac{1}{2}(n/9-41\nu n)>\alpha(G),$ where we used Fact~\ref{kisalfa} and that $1/9-41\nu>1/90.$ Both cases would imply the existence of a $K_4$ in $G,$
hence, we arrived at a contradiction.
\qedf

\medskip

Claim~\ref{bentifok} implies the following.

\begin{cor}\label{kor} Every vertex of $A$ can have up to $|A'|+\alpha(G)<16\nu n$ neighbors in $A,$ and similarly,
every vertex of $B$ can have up to $|B'|+\alpha(G)<16\nu n$ neighbors in $B.$ 
\end{cor}

But these degree bounds have even stronger consequences.

\begin{claim}\label{bentifok2}
If $v\in A$ then $deg(v, A)\le \alpha(G).$ Similarly, if $u\in B,$ then $deg(u, B)\le \alpha(G).$ 
\end{claim} 

\noindent {\bf Proof:}
We only sketch the proof as it is very similar to the proof of Claim~\ref{bentifok}. If some $v\in A$ has more than $\alpha(G)$ neighbors in $A,$ then one would
get a triangle inside $A.$ By Corollary~\ref{kor} we know that every vertex of $A,$ even those in $A',$ has almost $|B|/2$ neighbors in $B.$ Hence, in such a triangle we would have 
two vertices with more than $\alpha(G)$ common neighbors. But this would result in a $K_4,$ proving what was desired. The same argument works for $B$ as well.\qedf

\medskip

Another implication is the following.

\begin{claim}\label{meret}
There exists a non-negative integer $k\le 3\alpha(G)$ such that $n/2-k \le |A|, |B|\le n/2+k.$
\end{claim}

\noindent {\bf Proof:}
Let $u, v$ be two adjacent vertices from $A.$ By Claim~\ref{bentifok2} we have that $deg(u, A), deg(v, A) \le \alpha(G).$ Since
$\delta(G)\ge n/4,$ using Observation~\ref{fokos} this implies that $$|B|\ge deg(u, A)+ deg(v, A) -\alpha(G)\ge 2(n/4-\alpha(G))-\alpha(G)=n/2-3\alpha(G).$$ 
The same argument implies the lower bound $|A|\ge n/2-3\alpha(G).$ Since $A\cap B=\emptyset$ and $|A\cup B|=n,$ we proved what 
was desired.\qedf

\medskip

The above bounds for $|A|, |B|$ and the neighborhood structure of the vertices imply the following.

\begin{claim}\label{korok}
Neither $G[A],$ nor $G[B]$ can have a cycle of length 3, 5 or 7.
\end{claim}

\noindent {\bf Proof:} We will prove the statement for $G[A].$ Note that we have already proved that $G[A]$ cannot have a triangle in Claim~\ref{bentifok2}, so we only consider
the cases of $C_5$s and $C_7$s.

Assume first that $v_1, v_2, \ldots, v_5\in A$ such that $v_iv_{i+1}\in E$ for
$i=1, \ldots, 4.$ Set $H_i=N(v_i)\cap B$ for $i=1, \ldots, 5.$ Using that $G$ is $K_4$-free, we have that 
$H_1\cap H_2$ and $H_2\cap H_3$ both have at most $\alpha(G)$ vertices. Claim~\ref{bentifok2} implies, that
$|H_i|\ge n/4-\alpha(G)$ for $1\le i\le 5.$ 
By Claim~\ref{meret} we have $$|B-H_{i+1}| \le n/2+3\alpha(G)-(n/4-\alpha(G))= n/4 + 4\alpha(G).$$
Since $H_i\subseteq (B-H_{i+1})\cup (H_i\cap H_{i+1}),$ 
we obtain the upper bound $$|H_i|\le n/4+4\alpha(G) + \alpha(G) = n/4+5\alpha(G)$$ for $i=1, \ldots, 4.$
The above also imply that $$|B-(H_1\cup H_2)|\le n/2+3\alpha(G) -2(n/4-\alpha(G))+\alpha(G)=6\alpha(G).$$
Hence, $v_3$ can have at most $6\alpha(G)$ neighbors in $B-(H_1\cup H_2)$ and at most $\alpha(G)$ further neighbors in
$H_2,$ implying, that $|H_1\cap H_3|\ge n/4-8\alpha(G).$ 

The same way we obtain that $|H_3\cap H_5|\ge n/4-8\alpha(G).$ Since $|H_i|\le n/4+5\alpha(G),$ we have that 
$$|H_3-H_1|\le n/4+5\alpha(G)-(n/4-8\alpha(G))\le 13\alpha(G).$$ Therefore, $|H_5\cap (H_3-H_1)|\le 13 \alpha(G).$
These imply that $$|H_5\cap H_1|\ge |H_5\cap H_1\cap H_3|=|H_5\cap H_3|-|H_5\cap (H_3-H_1)|\ge n/4-21\alpha(G)>\alpha(G).$$ Hence, $G[A]$ must be $C_5$-free.

For proving the $C_7$-freeness of $G[A]$ assume that $v_1, v_2, \ldots, v_7\in A$ such that $v_iv_{i+1}\in E$ for $i=1, 2,\ldots, 6.$
Set $H_i=N(v_i)\cap B$ for every $i.$ Using the arguments of the previous case we have that $|H_3-H_1|\le 13\alpha(G),$  $|H_7\cap H_3|\ge n/4-21\alpha(G),$ and $|H_7 \cap (H_3-H_1)|\le 13\alpha(G).$ These imply that 
$$|H_7\cap H_1|\ge |H_7\cap H_3\cap H_1|=|H_7\cap H_3|-|H_7\cap (H_3-H_1)|\ge n/4-34\alpha(G)>\alpha(G),$$ that is,
$G[A]$ must be $C_7$-free. \qedf

\medskip

The following is Lemma 7.1 in~\cite{LR}, we omit the proof.

\begin{lemma}\label{alfan}
Every graph $F$ not containing a cycle of length 3, 5, or 7 satisfies the inequality
$$e(F) \le \alpha(F)^2.$$
\end{lemma} 

\medskip

\noindent {\bf Step 5.} We divide $A$ and $B$ into disjoint subsets as follows: $A=I_A\cup L_A \cup M_A$ and $B=I_B\cup L_B\cup M_B.$ 
Here $I_A$ includes all those vertices $v\in A$ for which $deg(v, B)> |B|/2+ 4\alpha(G),$ $L_A$ includes those vertices
$v\in (A-I_A)$ for which $(|B| + \alpha(G))/2<deg(v, B)\le  |B|/2+ 4\alpha(G).$ Finally, $M_A=A-(I_A\cup L_A).$ 
The subsets $I_B, L_B$ and $M_B$ of $B$ are defined analogously, only the roles of $B$ and $A$ are interchanged. 

\begin{claim}\label{fgetlen}
The subset $I_A\cup L_A$ is an independent set in $G[A],$ and similarly, $I_B\cup L_B$ is an independent set in $G[B].$ Moreover,
$I_A$-vertices have no neighbor in $A$ and $I_B$-vertices have no neighbor in $B.$  
\end{claim}

\noindent {\bf Proof:} Let $u, v\in I_A\cup L_A$ any two distinct vertices. By definition their degree sum is greater
then $|B|+\alpha(G),$ hence, $uv\not\in E$ by the $K_4$-freeness of $G.$ A similar statement holds for any two distinct vertices
of $I_B\cup L_B.$ 

Using that no vertex of $A$ has more than $\alpha(G)$ neighbors in $A,$ and $|B|\le n/2+3\alpha(G)$ by Claim~\ref{meret},
every $v \in A$ has at least $n/4-\alpha(G)>|B|/2-3\alpha(G)$ neighbors in $B.$ Hence, $u\in I_A$ cannot have any
neighbor in $A$ by Observation~\ref{fokos}. The same argument works for $I_B,$ too.
\qedf

\medskip

Next we perform an operation which may increase the number of edges without creating any $K_4$ in the new graph, 
which we denote by $G'$ (such an operation is used in~\cite{FLZ}).
The details are as follows.

\begin{itemize}

\item For every $v\in L_A$ delete the edges that connect $v$ with any other vertex in $A,$ and include all edges
of the form $vu$ where $u\in B.$ We do the same for the vertices of $L_B$: delete the edges that go inside, and
include all edges that go in between $A$ and $B.$

\item For every $v\in I_A$ include all edges $uv$ with $u\in B.$ Similarly, include all edges $uv$ where $u\in I_B$ and
$v\in A.$

\end{itemize}

Let us remark, that as $G[A], G[B]$ had no cycles of length 3, 5 or 7, this also applies for $G'[A]$ and $G'[B],$
since we could only delete edges from inside $A$ and $B.$

\medskip

\begin{claim}\label{nincsK4}
The new graph $G'$ cannot have less edges than $G,$ moreover, $G'$ is $K_4$-free.
\end{claim}

\noindent {\bf Proof:} The vertices of $L_A$ had at most $|B|/2+4\alpha(G)$ neighbors in $B$ by definition, and at most $\alpha(G)$
neighbors in $A,$ all those belong to $M_A.$ 

Hence, $e(G[M_A, L_A])\le \alpha(G) |L_A|.$ We can get the bound $e(G[M_B, L_B])\le \alpha(G) |L_B|$ analogously. That is, during
the operation we may lose at most $\alpha(G) |L_A|$ edges inside $A,$ and similarly, at most $\alpha(G) |L_B|$ edges inside $B.$

Next we estimate the number of new edges. 
Let $v\in L_A$ be arbitrary. It had at most $(|B|/2+4\alpha(G))$ neighbors in $M_B$ by definition, hence, we included at least 
$|M_B|-(|B|/2+4\alpha(G))\ge |B|-\alpha(G)-(|B|/2+4\alpha(G))=|B|/2-5\alpha(G)$ new edges at $v.$ 
Analogously, for every $v\in L_B$ we included at least $|A|/2-5\alpha(G)$ new edges. 
 
Hence, the total change in the number of edges is at least $$(|B|/2-5\alpha(G))\cdot |L_A|+(|A|/2-5\alpha(G))\cdot |L_B|- (|L_A|+|L_B|)\cdot \alpha(G).$$
This expression is easily seen to be non-negative, since $\alpha(G)\ll n.$
  
Suppose on the contrary now, that $G'$ has a $K_4.$ First notice that such a $K_4$ cannot contain two vertices
from $I_A\cup L_A$ or from $I_B\cup L_B,$ since these are independent sets. We cannot have a $K_4$ containing 
one vertex from $I_A\cup L_A$ and another vertex from $I_B\cup L_B,$ since in $G'$ such a vertex does only have
neighbors from the opposite part. By the triangle-freeness of $G[A]$ and $G[B]$ the subgraphs $G'[A]$ and $G'[B]$
must also be triangle-free, hence, no $K_4$ can contain any vertex from $I_A\cup L_A\cup I_B\cup L_B.$ 
There is only one possibility left: when all the four vertices belong to $M_A\cup M_B.$ But such a $K_4$ would be
a $K_4$ in $G$ as well, so we arrived at a contradiction.
\qedf

\medskip

\begin{claim}\label{fuggmeret}
We have $e(G'[M_A])\le (\alpha(G)-(|I_A|+|L_A|))^2$ and  $e(G'[M_B])\le (\alpha(G)-(|I_B|+|L_B|))^2$
\end{claim}

\noindent {\bf Proof:} The statement follows from Lemma~\ref{alfan} and the facts that $I_A\cup L_A$ and $I_B\cup L_B$ are independent
sets in $G'[A],$ respectively, $G'[B],$ and that $e(G'[M_A, I_A\cup L_A])=e(G'[M_B, I_B\cup L_B])=0.$ Hence, the
largest independent set in $G'[M_A]$ can have at most $\alpha(G)-|I_A|-|L_A|$ vertices, and a similar 
statement holds for
the largest independent set of $G'[M_B].$ Using Lemma~\ref {alfan} we obtain the claimed bounds. 
\qedf

\medskip

We are ready to finish the proof of Theorem~\ref{fottetel}. Putting together all the above, we can give the desired upper
bound for the number of edges in $G'.$ As it turns out, the calculation is somewhat simpler if we bound $2e(G').$ Set $a=|A|,$ $b=|B|,$ $m_A=|M_A|$
and $m_B=|M_B|.$ Without loss of generality we assume that $|A|=n/2+k$ and $|B|=n/2-k$ for some $0\le k\le 3\alpha(G),$
given by Claim~\ref{meret}.
We define two functions: 
$$f_A(x)= x \frac{b+\alpha(G)}{2} + (a-x)b + 2(\alpha(G) - (a-x))^2$$ and 
$$f_B(x)= x \frac{a+\alpha(G)}{2} + (b-x)a+2(\alpha(G) - (b-x))^2.$$  Observe, that $f_A(m_A)-(\alpha(G) - (a-m_A))^2$ is an upper bound for the number of edges incident 
to vertices of $A$ and  $f_B(m_B)-(\alpha(G) - (b-m_B))^2$ is an upper bound for the number of edges incident 
to vertices of $B.$ Hence, $2e(G')\le f_A(m_A)+f_B(m_B).$ 

Elementary calculus shows that $f_A(x)$ and $f_B(x)$ are monotone decreasing, and reach their maximum at $m_A=a-\alpha(G)=n/2+k-\alpha(G)$ and $m_B=b-\alpha(G)=n/2-k-\alpha(G),$ respectively. 
Plugging in these values we get the following upper bounds:
$$f_A(n/2+k-\alpha(G))= (n/2 +k -\alpha(G))\frac{n/2-k+\alpha(G)}{2}+\alpha(G) (n/2-k)$$ and
$$f_B(n/2-k-\alpha(G))= (n/2 -k -\alpha(G))\frac{n/2+k+\alpha(G)}{2}+\alpha(G) (n/2+k).$$

Simple calculation shows, that $$f_A(n/2+k-\alpha(G))=f_B(n/2-k-\alpha(G))=\frac{n^2}{8}+\frac{\alpha(G)n-\alpha(G)^2}{2}-\frac{k^2}{2},$$ 
therefore $$e(G)\le \frac{n^2}{8}+\frac{\alpha(G) n}{2}-\frac{\alpha(G)^2}{2} -\frac{k^2}{2}\le  \frac{n^2}{8}+\frac{\alpha(G) n}{2}-\frac{\alpha(G)^2}{2}.$$
This finishes the proof of the theorem. 
\qedf

\end{document}